\documentclass[11pt]{article}
\usepackage{xy}
\usepackage{amsmath,amsthm}
\usepackage{amssymb,latexsym}
\usepackage{mathrsfs}
\usepackage{makeidx}
\makeindex

\def\lim{\mathop{\rm lim}}

\def\phi{\varphi}

\def\pf{{\it Proof:}~}
\def\qed{$\,_{\square}\\$}
\newtheorem{theorem}{Theorem}
\newtheorem{lemma}[theorem]{Lemma}
\newtheorem{proposition}[theorem]{Proposition}
\newtheorem{definition}[theorem]{Definition}
\newtheorem{corollary}[theorem]{Corollary}
\newtheorem{problem}[theorem]{Problem}
\newtheorem{exercises}[theorem]{Exercises}
\newtheorem{exercise}[theorem]{Exercise}
\newtheorem{remark}[theorem]{Remark}
\newtheorem{claim}[theorem]{Claim}
\newtheorem{examples}[theorem]{Examples}
\newtheorem{question}[theorem]{Question}
\newtheorem{conjecture}[theorem]{Conjecture}

\newcommand{\begintheorem}{\begin{theorem}}
\newcommand{\beginlemma}{\begin{lemma}}
\newcommand{\beginexercise}{\begin{exercise}}
\newcommand{\beginexercises}{\begin{exercises}}
\newcommand{\beginproposition}{\begin{proposition}}
\newcommand{\begindefinition}{\begin{definition}}
\newcommand{\begincorollary}{\begin{corollary}}
\newcommand{\beginproblem}{\begin{problem}}
\newcommand{\beginremark}{\begin{remark}}
\newcommand{\beginclaim}{\begin{claim}}
\newcommand{\beginassumptions}{\begin{assumptions}}
\newcommand{\beginexamples}{\begin{examples}}
\newcommand{\beginquestion}{\begin{question}}
\newcommand{\beginsassumptions}{\begin{sassumptions}}
\newcommand{\beginsassumption}{\begin{sassumption}}
\newcommand{\beginconjecture}{\begin{conjecture}}

\author{Robert Gulliver and Guoyi Xu \\School  of Mathematics\\
University of Minnesota\\ Minneapolis, MN 55455}

\title{\huge \bf Examples of hypersurfaces flowing by curvature
in a Riemannian manifold}
\begin{document}
\maketitle
\begin{abstract}
This paper gives some examples of hypersurfaces $\phi_t(M^n)$
evolving in time with speed determined by functions of the normal
curvatures in an $(n+1)$-dimensional hyperbolic manifold; we
emphasize the case of flow by harmonic mean curvature. The
examples
converge to a totally geodesic submanifold of any dimension from
$1$ to $n$, and include cases which exist for infinite time.
Convergence to a point was studied by Andrews, and only
occurs in finite time.  For dimension $n=2,$ the destiny of
any harmonic mean curvature flow is strongly influenced by the
genus of the surface $M^2$.\\[3mm]
Mathematics Subject Classification: 35K15, 53C44
\end{abstract}

%
%
\section{Background}
Unless otherwise mentioned, all Riemannian manifolds in this
article are connected and complete. Let $M^n$ be a smooth,
connected, orientable compact manifold of dimension $n\geq 2$,
without
boundary, and let $(N^{n+1},g^N)$ be a smooth connected Riemannian
manifold.  $\sigma^N$ is any sectional curvature of $N^{n+1}$,
$\mathscr R$ is the Riemann tensor of $N^{n+1}$, and $\nabla^N$ is
the
Levi-Civita connection corresponding to $g^N$.  For a hyperbolic
manifold, $\sigma^N \equiv -1$.  When an index such
as $i$ is repeated in one term of an expression, summation
$1\leq i \leq n$ is indicated.

Suppose $\varphi_0:M^n\rightarrow N^{n+1}$ is a smooth immersion
of an oriented manifold $M^n$ into $N^{n+1}$; write $\vec{v}$
for the induced normal vector to $\varphi_0(M)$. The second
fundamental form of $M$ is a covariant tensor which we represent
at each point by a matrix $A$, where the entry
$A_{ij}=h_{ij}=\Big<\nabla_{\frac{\partial}{\partial
x_i}}^{N}\vec{v},\frac{\partial}{\partial x_j}\Big>_{g^N}$. The
Weingarten tensor is given by the matrix $\mathscr{W}$,
whose entry $\omega_{i}^{k}=h_{ij}g^{jk}$ and
$\{g^{jk}\}$
is the pointwise inverse matrix of
$\{g_{jk}\}$.\\
We seek a solution $\varphi:M^n\times[0,T)\rightarrow N^{n+1}$ to
an equation

\begin{equation}\label{eq:1}
\frac{\partial}{\partial t}\varphi(x,t)=
-f(\lambda(\mathscr{W}(x,t)))\vec{v}(x,t)\\
\end{equation}
\[\varphi(x,0)=\varphi_0 (x)\]

\noindent
where $F(x,t)=f(\lambda(\mathscr{W}(x,t)))$ and $f$ is a smooth
symmetric function, where $\vec{v}(x,t)$ is the outward normal
vector to $\phi(M^n,t)$.  $\mathscr W(x,t)$ is the Weingarten
matrix of $\varphi(M^n,t)$ in $N^{n+1}$, and
$\lambda(\mathscr{W})$ is the set of eigenvalues
$(\lambda_1,\ldots,\lambda_n)$ of $\mathscr{W}$. Define
$\phi_t(x)=\phi(x,t)$, then $(\lambda_1,\ldots,\lambda_n)$ are the
principal curvatures of the hypersurface
$M_t\stackrel{\vartriangle}{=}\varphi_t(M)\subset N$. \\ 

For example,
(\ref{eq:1}) becomes Mean Curvature Flow when
$f(\lambda)=\sum_{i}\lambda_i$ (see \cite{CRM}, \cite{LMCF}). 

Consider the solution $\phi:M^n\times[0,T)\rightarrow N^{n+1}$ of
the following equations:

\begin{equation}\label{HMF}
\frac{\partial}{\partial t}\varphi(x,t)=
-\Big(\sum_i\lambda_i^{-1}\Big)^{-1}\vec{v}(x,t)\\
\end{equation}
\[\varphi(x,0)=\varphi_0 (x)\]
Such a solution $\varphi(x,t)$ is Harmonic Mean Curvature Flow;
$f(\lambda)=(\sum_i \lambda_i^{-1})^{-1}$ is the harmonic mean of
the numbers $\lambda_1,\dots,\lambda_n$.

B. Andrews proved the following theorem in $\cite{AR}$:

%
\begintheorem \label{Andrews}
{Let $M^n$ and $\phi_0$ be assumed as at the beginning of this
paper, and that the Riemannian manifold $(N^{n+1},g^{N})$ 
satisfies the following conditions:
\[-K_1\leq\sigma^{N}\leq K_2,\quad
\left|\nabla^{N}R^N\right|_{g^N}\leq L\]
for some nonnegative constants $K_1$, $K_2$ and $L$.\\
Assume every principal curvature $\lambda_i$ of $\phi_0$
satisfies the following condition:
\[\lambda_i>\sqrt{K_1}\]
Then there exists a unique smooth solution to (\ref{HMF}) on a
maximal time interval $[0,T)$, $T<\infty$, and the immersion
$\phi_t$ converges uniformly to a round point $p$ in $N^{n+1}$ as
$t$ approaches $T$.}
\end{theorem}

Also, we have the following theorem, to appear in $\cite{GYX}$:

%
\begintheorem \label{main-2}
{Let $M^n$ be a smooth, connected, orientable compact manifold of
dimension
$n\geq 2$, without boundary. Assume $N^{n+1}$ is a non-positively
curved, simply-connected smooth manifold, and suppose
$\phi_0: M^n\rightarrow N^{n+1}$ is a smooth immersion of $M^n$.
Assume every principal curvature of $\phi_0(M)$ is positive.
Then there exists a unique smooth solution to (\ref{HMF}) on a
maximal time interval $[0,T)$, $T<\infty$, and the immersion
$\phi_t$ converges uniformly to a round point $p$ in $N^{n+1}$ as
$t$ approaches $T$.}
\end{theorem}

In the rest of this paper, except for Section 6, and unless
otherwise mentioned, we consider harmonic mean curvature flow and
let $f(\lambda)=(\sum_i \lambda_i^{-1})^{-1}$.
We provide two specific examples of harmonic mean curvature flow
for infinite time: 
in section 2, with dimension reduction in the limit, and in
section 3, with the limit manifold of the same dimension as $M$.
Note these examples in section 2 and section 3 provide barriers
for harmonic mean curvature flow in Riemannian manifolds; further
applications will be addressed in \cite{GYX}. We discuss the limit
behavor of the harmonic mean curvature flow at infinite time in
section 4. Then we treat the special consequences of the
Gauss-Bonnet theorem for $2$-dimensional surfaces in section 5,
and turn to examples of more general flows by functions of normal
curvatures in section 6.

We would like to thank Gerhard Huisken for interesting
discussions, and in particular for the observation that there are
no examples in the literature for convergence of a compact
hypersurface flowing by harmonic mean curvature in infinite time
to a set of positive dimension.  And we also would like to thank
the referee for pointing out a gap in the earlier version of this
paper.

%
%
\section{The dimension-reduction example}

In this section, we give an example where $\phi_t$ converges to
$\phi_{\infty}$ in the $C^{\infty}$ topology but the dimension of
$M_{\infty}=\phi_{\infty}(M)$ is less than the dimension of
$M_t$; i.e. there is {\bf dimension reduction}.

%
\begintheorem\label{dimred}
{
Let $N^3$ be a hyperbolic manifold containing an embedded 
closed geodesic $M_\infty$. Then there is a flow 
$\phi_t: M^2\rightarrow N^3$ by harmonic mean curvature, where 
$M^2$ is a
torus, which converges to $M_{\infty}$ as $t\rightarrow +\infty$.}
The flow consists of immersions $\phi_t$, which become embedded
for $t$ sufficiently large.
\end{theorem}

For example, we may let the ambient manifold $N$ be 
$H^3/\mathbb {Z}$, where $H^3$
is hyperbolic space, represented as the Poincar\'e half space
$\mathbb ({R}^3)^{+}=\{(x,y,z)|(x,y,z)\in \mathbb R^3, z>0\}$
with the metric $g^N_{ij}=\frac{1}{z^2}\delta_{ij}$ ($\delta_{ij}= \delta_i^j=$ Kronecker delta), and the
$\mathbb Z$ action $f: \mathbb {Z} \times H^3 \to H^3 $ is defined
as:
\[f(k)(x,y,z)=2^{k}(x,y,z).\]
Recall that $f(k)$ is an isometry of $H^3$ for each 
$k\in \mathbb Z$.

Now we let $N$ be the quotient manifold of $H^3$ under the
$\mathbb Z$-action, with fundamental domain
$\{(x,y,z)|\, 1\leq \sqrt{x^2+y^2+z^2}\leq 2\}$. Then $M_\infty=$
the positive $z$-axis, modulo $f(1)$, is a closed geodesic in $N$.\\

\pf
Let $\psi_0:\mathbb S^1 \to N$ be an embedding as the given closed
geodesic curve $M_\infty$ in $N$.  We choose a 
unit vector field $w(x)$ in $(T_{x}\psi_0)^{\perp}$. Then for
$r>0$, we define
\[\psi(x,\theta,r)=\psi_r(x,\theta):\mathbb S^1\times\mathbb
S^1\to N^3\]
by
\[\psi(x,\theta,r)=\psi_r(x,\theta)=\gamma(x,\theta,r),\]
where $\gamma(x,\theta,\cdot)$ is the unit-speed geodesic in $N$
with $\gamma(x,\theta,0)= \psi_0 (x)$ and
$\frac{d}{dr}\gamma(x,\theta,r)=\vec{N}(x, \theta)$ at $r=0$.
Here $\vec{N}(x,\theta)$ is the unit tangent vector in
$T_{\psi_0(x)}N^3$
which is perpendicular to $T_x\psi_0$ and makes the angle
$\theta$ with $w(x).$  Then
$\psi_r(\mathbb S^1\times \mathbb S^1)$
has two principal curvatures:
\[\lambda_1(r) \equiv \tanh r, \quad \lambda_2(r) \equiv \coth
r.\]
In fact, for $i=1,2$, $\lambda_i(r)$ is the logarithmic
derivative of the length of a Jacobi field, and hence satisfies
the Ricatti equation $\lambda_i'(r)+(\lambda_i(r))^2=1$.

We have constructed a one-parameter family of immersions
$\psi_r : M \to N$, $-\infty < r < \infty,$ with two principal
curvatures:
$\lambda_1(r) \equiv \tanh r$ and $\lambda_2(r) \equiv \coth r$.
It may be observed that $\psi_r$ is an embedding for $r$
sufficiently small.\\

Now consider the harmonic mean curvature flow
$\phi_t=\psi_{r(t)}: M \to N$, with initial conditions
$\phi_0 = \psi_{r_0}$, $r(0)= r_0$, where $r_0$ is some fixed
positive constant.  The speed must satisfy:
\[\frac{\partial r}{\partial t}=
\Big<\frac{\partial\gamma}{\partial r}\frac{\partial r}{\partial t},
\vec{v}\Big> =
\Big<\frac{\partial{\gamma(x,r)}}{\partial t},\vec{v}\Big> \]
\[=\Big<\frac{\partial{\psi(x,r)}}{\partial t},\vec{v}\Big> =
\Big<\frac{\partial{\phi(x,t)}}{\partial t},\vec{v}\Big>\]
\[=\Big<-F\vec{v},\vec{v}\Big> =-F(\lambda_1, \lambda_2)\]
\[=-\frac{1}{\lambda_1^{-1}+ \lambda_2^{-1}}=
-\frac{\sinh r\cosh r}{(\sinh r)^2 +(\cosh r)^2}.\]
In the first equation we use the fact  
$\frac{\partial \gamma}{\partial r}= \vec{v}$;
in the third equation we use the definition of $\psi_r$, where
$\vec{v}=\vec{N}(x,\theta)$ is the outward normal vector of 
$\psi_r(M)$ at $(x,\theta)\in \mathbb S^1\times\mathbb S^1$.

Solving, we find
$$r(t) = \frac{1}{2} \sinh^{-1} \Big(e^{-t} \sinh 2r_0 \Big).$$
Note that $r(t) \to 0$ as $t \to \infty$.  \qed

%
%
\section{The no-dimension-reduction example} 

In this section, we give an example in which $M_t$ converges to
$M_{\infty}$ in the $C^{\infty}$ topology and the dimension of
$M_{\infty}$ is the same as the dimension of $M_t$, i.e. there is
{\bf no dimension reduction.}

%
\begintheorem\label{nodimred}
{There is a compact surface $M^2$ of genus $2$, a hyperbolic
manifold $N^3$ diffeomorphic to $M\times \mathbb R$, a totally
geodesic embedding $\psi_0: M\rightarrow N$ and a flow by harmonic
mean curvature $\phi_t: M\rightarrow N$ such that as
$t\rightarrow +\infty$, $\phi_t(M) \rightarrow \psi_0(M)$
smoothly.}
\end{theorem}

\pf
{ Let $\Omega$ be a regular geodesic octagon in the hyperbolic
plane $H^2$, with angles $\pi/2$,
and thus area $4\pi$.  Label the edges as
$$\beta_1, \alpha_1', -\beta_1', -\alpha_1, \beta_2, \alpha_2',
-\beta_2', -\alpha_2,$$ 
in that order, where the signs indicate
orientation.  Let $A_1$ be the orientation-preserving isometry of
$H^2$ which maps the oriented geodesic segments $\alpha_1$ to
$\alpha_1'$;  $A_2$ maps $\alpha_2$ to $\alpha_2'$; $B_1$ maps
$\beta_1$ to $\beta_1'$; and $B_2$ maps $\beta_2$ to $\beta_2'$.
The group $G$ of isometries of $H^2$ generated by $A_1, A_2$ and
$B_1$ also includes $B_2$.  $G$ is isomorphic to the fundamental
group of the compact surface of genus $2$. (See pp. 95--98 in Katok
\cite{FG} for the arithmetic properties of the group $G$.)

Let $\psi_0:H^2 \to H^3$ be an embedding as a totally geodesic
surface in $H^3$.  The isometries in $G$ extend in a well-known
fashion to isometries of $H^3$, leaving the distance from
$\psi_0(H^2)$ invariant.

Choose a unit normal vector field $\vec{N}$ to
$\psi_0(H^2)$.  Define $\psi(\cdot,r): H^2 \to H^3$ by
$\psi(x,r)=\psi_r(x)=\gamma(x,r)$ and $\psi(x,0)=\psi_0(x)$,
where $\gamma(x,\cdot)$ is the unit-speed geodesic
in $H^3$ with $\gamma(x,0)=x$ and
$\frac{\partial}{\partial r}\gamma (x,0)=\vec{N}(x)$.

Then $\psi_r(H^2)$ is totally umbilic, with normal curvatures
$\lambda(r) \equiv \tanh r$.  In fact, $\lambda(r)$
satisfies the Ricatti equation $\lambda'(r)+(\lambda(r))^2=1$,
with the initial condition $\lambda(0)=0$.

Now let the group $G$ act by isometries on $H^2$ and on $H^3$.
The quotient $H^2/G = M^2$ is a compact surface of genus $2$, 
with fundamental domain $\Omega$, and
the quotient $H^3/G = N^3$ is a noncompact hyperbolic manifold
diffeomorphic to $M \times {\mathbb{R}}$.  The group $G$ acting on
$N$ preserves each of the hypersurfaces $\psi_r(H^2)$.  We have
constructed a one-parameter family of totally umbilic embeddings
$\psi_r : M \to N$, $-\infty < r < \infty,$ with normal curvatures
$\equiv \tanh r$.

Now consider the harmonic mean curvature flow $\phi_t: M \to N$,
with initial conditions $\phi_0 = \psi_{r_0}$, where $r_0$ is some
fixed positive constant.  The speed must satisfy

\[\frac{\partial r}{\partial t}=
\Big<\frac{\partial \gamma}{\partial r}\frac{\partial r}{\partial
t}, \vec{v}\Big> =
\Big<\frac{\partial{\gamma(x,r)}}{\partial t},\vec{v}\Big> =
\Big<\frac{\partial{\psi(x,r)}}{\partial t},\vec{v}\Big>  \]
\[= \Big<\frac{\partial{\phi(x,t)}}{\partial t},\vec{v}\Big> =
\Big<-F\, \vec{v},\vec{v}\Big> =-F(\lambda_1, \lambda_2)\]
\[=-\frac{1}{\lambda_1^{-1} + \lambda_2^{-1}} =-\frac12 \tanh r.\]

\noindent
In the first equation we use the fact
$\frac{\partial \gamma}{\partial r}=\vec{N}(x)=\vec{v}$.
In the third equation we use the definition of $\psi_r$, where
$\vec{v}$ is the outward normal vector of $\psi_r$.

Solving, we find
$$ r(t) = \sinh^{-1} \Big(e^{-t/2} \sinh r_0 \Big).$$
Note that $r(t) \to 0$ as $t \to \infty$.
} \qed

%
%
\section{The limit behavior of harmonic mean curvature flow at
infinite time}

In this section, we will give a sufficient condition where 
harmonic mean curvature flow will exist forever, and discuss the
limit behavior.  Let $\varphi_t:M\to N$ be an immersion of $M^n$
into a hyperbolic manifold $N^{n+1}.$

\begindefinition
{We define the following notation:
\[\dot{F}^{kl}=\frac{\partial F}{\partial h_{kl}},\quad \ddot{F}^{kl,pq}= \frac{\partial^2 F}{\partial h_{kl} \partial h_{pq}}, \quad \dot{H}_k^i= \frac{\partial H}{\partial \omega_i^k},\quad \ddot{H}_{r,k}^{s,i}= \frac{\partial^2 H}{\partial \omega_i^k\partial \omega_s^r},\quad \mathscr
R_{ij}=\mathscr R_{i0j0},\]
where $0$ appearing as a tensor index represents the
normal vector $\vec{v}$ of $\varphi(M)$ in $N$. For
any $W:M\to \mathbb R$, we define:

\[\mathscr L(W)=\dot{F}^{kl}\nabla_k\nabla_lW.\]
}
\end{definition}

Recall from Andrews \cite{AR} that $\mathscr L$ is elliptic as
long as $\phi_t(M)$ remains locally strictly convex. 

%
\begintheorem \label{hyperbolic infinite time}
{If $N^{n+1}$ is a hyperbolic manifold, $F(x)
< \frac{1}{n}$ for any $x\in M$,
then $\phi_t(M)$ remains locally convex and $F(x,t)< \frac{1}{n}$
for any $x\in M$, $t\in [0, +\infty)$,
$\lim_{t\rightarrow \infty} F(x,t)= 0$, and the harmonic mean
curvature flow exists for all $t$ in $[0,+\infty)$.}
\end{theorem}

\pf
{ By Andrews \cite{AR}, using a curvature coordinate system at one
point, we have the following formula:\\
\[\frac{\partial F}{\partial t}=
\mathscr L(F)+F<\dot F,(\mathscr W^2)>+ F<\dot F^{ij},(\mathscr
R_{ij})>\]

\begin{equation} \label{formula}
{=\mathscr L(F)+ \sum_{i}F \,\frac{\partial f}{\partial
\lambda_i} \,(\lambda_i^2+\mathscr R_{ii})
}
\end{equation}

\[\leq\mathscr L(F)+ F^3(n -\sum_i {\lambda_i^{-2}})\]
\[\leq\mathscr L(F)+ F^3\Big(n -\frac{1}{n}F^{-2} \Big).\]

Consider the ODE
\[\frac{\partial \tilde{F}}{\partial t}=
\tilde{F}^3(n-\frac{1}{n}\tilde{F}^{-2}),\]
\[\tilde{F}(0)=\max_{x\in M}F(x,0).\]

Solving the above ODE, we get
$\widetilde F(t)^{-2}-n^2=(\widetilde F(0)^{-2}-n^2)e^{2t/n}.$
Because $0< \widetilde F(0)= \max_{x\in M^n} F(x,0)< \frac{1}{n}$,
we get $\lim_{t\rightarrow \infty} \tilde{F}(t)= 0$.

By the maximum principle, $F(x,t)\leq \tilde{F}(t)< \frac{1}{n}$,
for all $x\in M$, $t\in [0, +\infty)$, and therefore
$\lim_{t\rightarrow \infty} F(x,t)= 0$.

On the other hand, we have the following estimate by the above
evolution equation of $F$:
\[\frac{\partial F}{\partial t}\geq
\mathscr L(F)+ F^3(-\sum_i {\lambda_i^{-2}})\geq
\mathscr L(F)- F.\]

Now consider the ODE
\[\frac{\partial \widehat{F}}{\partial t}= -\widehat{F},\]
\[\widehat{F}(0)=\min_{x\in M}F(x,0).\]

\noindent
Then by the maximum principle again, we get for all $x\in M, t\in
[0, +\infty)$: \[ F(x,t)\geq \widehat{F}(t) =
\min_{x\in M}F(x,0)\,e^{-t}>0\]

In particular, $\phi_t(M)$ remains convex for all $t$.

Finally, we have the following estimate of $H$. By Andrews'$\cite{AR}$:
\[\frac{\partial}{\partial t}\omega_i^{r}=
\dot{F}^{kl}\nabla_k\nabla_l\omega_i^r+
\ddot{F}^{kl,pq}(\nabla_i h_{kl})(\nabla_j h_{pq})g^{jr}\]
\[+\dot{F}^{kl}(h_{ml}\omega_k^m)\omega_i^r + 
\dot{F}^{st}\mathscr R_{st}h_{ij}g^{jr} + 
2\dot{F}^{pm}g^{tr}\omega_m^q\mathscr R_{piqt}\]
\[-\dot{F}^{pq} (g^{tr}\omega_i^s\mathscr R_{psqt}+ 
g^{ts}\omega_s^r\mathscr R_{piqt}) + 
\dot{F}^{pq}g^{tr}(\nabla_i\mathscr R_{tpq0}-\nabla_p\mathscr R_{qit0})\]

Now referring to the last five terms above, we define: 
\[(I)=\dot{H}_r^i\dot{F}^{kl}(h_{ml}\omega_k^m)\omega_i^r\ ,\  
(II)=\dot{H}_r^i\dot{F}^{st}\mathscr R_{st}h_{ij}g^{jr}\]
\[(III)=2\dot{H}_r^i\dot{F}^{pm}g^{tr}\omega_m^q\mathscr R_{piqt} \ , 
\ (IV)=-\dot{H}_r^i(\dot{F}^{pq}g^{tr}\omega_i^s\mathscr R_{psqt}+
\dot{F}^{pq}g^{ts}\omega_s^r\mathscr R_{piqt})\]
\[(V)=\dot{H}_r^i\dot{F}^{pq}g^{tr}
(\nabla_i\mathscr R_{tpq0}-\nabla_p\mathscr R_{qit0})\]

\noindent
then 
\[\frac{\partial}{\partial t}H=
\dot{H}_r^i(\frac{\partial}{\partial t}\omega_i^r)\]
\[=\dot{H}_r^i(\dot{F}^{kl}\nabla_k\nabla_l\,\omega_i^r)+
\dot{H}_r^i\ddot{F}^{kl,pq}(\nabla_i h_{kl})(\nabla_j h_{pq})g^{jr}+
(I)+\cdots+(V)\]

Note 
\[\dot{F}^{kl}\nabla_k\nabla_lH=
\dot{F}^{kl}\nabla_k(\dot{H}_r^i\nabla_l\omega_i^r)
=\dot{F}^{kl}\ddot{H}_{r,\tilde{r}}^{i,\tilde{i}}
(\nabla_k\omega_{\tilde{i}}^{\tilde{r}})
(\nabla_l\omega_i^r)+\dot{F}^{kl}\dot{H}_r^i\nabla_k\nabla_l\omega_i^r.\]

Define \[(J)=\dot{H}_r^i\ddot{F}^{kl,pq}(\nabla_i h_{kl})
(\nabla_j h_{pq})g^{jr}-
\dot{F}^{kl}\ddot{H}_{r,\tilde{r}}^{i,\tilde{i}}
(\nabla_k\omega_{\tilde{i}}^{\tilde{r}})(\nabla_l\omega_i^r);\]
we get
\[\frac{\partial}{\partial t}H=\mathscr L(H)+(J)+(I)+\cdots+(V).\]

It is straightforward to get
\[(I)+(II)=H[<\dot{F},(\mathscr W^2)>+\dot{F}^{ij}
\mathscr R_{i0j0}]\leq nF^2H\leq \frac{1}{n}H\]
\noindent
and 
\[(V)= \frac{\partial f}{\partial \lambda_i}
(\nabla_j\mathscr R_{jii0}-\nabla_i\mathscr R_{ijj0})= 0.\]

Choose a curvature coordinate system around one point; then we 
could do the following 
calculation:
\[(J)=\ddot{F}^{kl,pq}(\nabla_ih_{kl})(\nabla_ih_{pq})\]
But by the Lemma 2.22 in \cite{AE}, we know $F$ is concave from 
the fact that $f$ is concave. So we get $(J)\leq 0$. 

Now 
\[(III)+(IV)=2\dot{H}_r^i\dot{F}^{pm}g^{tr}\omega_m^q\mathscr R_{piqt}-
\dot{H}_r^i(\dot{F}^{pq}g^{tr}\omega_i^s\mathscr R_{psqt}+
\dot{F}^{pq}g^{ts}\omega_s^r\mathscr R_{piqt})\]
\[=2\delta_r^i\frac{\partial f}{\partial \lambda_p}
\delta_p^m\delta_t^r\lambda_q\delta_q^m\mathscr R_{piqt}- 
\delta_r^i \big (\frac{\partial f}{\partial \lambda_p}
\delta_p^q\delta_t^r\lambda_i\delta_i^s\mathscr R_{psqt}+ 
\frac{\partial f}{\partial \lambda_p}\delta_p^q\delta_t^s\lambda_s\delta_s^r
\mathscr R_{piqt} \big)\]
\[=2\mathscr R_{prpr}\frac{\partial f}{\partial \lambda_p}(\lambda_p-\lambda_r)
=2\sum_{p<r}\mathscr R_{prpr} \big(\frac{\partial f}{\partial \lambda_p}-
\frac{\partial f}{\partial \lambda_r}\big) (\lambda_p-\lambda_r)\]
\[=\big(\sum_k\lambda_k^{-1}\big)^{-2}\cdot\sum_{i,j}(-\mathscr R_{ijij})
\cdot(\lambda_i-\lambda_j)^2(\lambda_i+
\lambda_j)\cdot\lambda_i^{-2}\lambda_j^{-2}\]
\[\leq \sum_{i,j}(\lambda_i+\lambda_j)\cdot
\Big(\frac{\lambda_i^{-1}-\lambda_j^{-1}}{\sum_k \lambda_k^{-1}}\Big)^2 
\leq \sum_{i,j}(\lambda_i+\lambda_j)= 2nH\]

We have the following inequality for $H$ by the above estimates:
\[\frac{\partial H}{\partial t}\leq \mathscr L(H)+ 
\Big(2n+\frac{1}{n}\Big)H.\]

Now consider the ODE
\[\frac{\partial \widehat{H}}{\partial t}= 
\Big(2n+ \frac{1}{n}\Big)\widehat{H},\]
\[\widehat{H}(0)=\max_{x\in M}H(x,0).\]

\noindent
Then by the maximum principle again, we get for all 
$x\in M, t\in [0, +\infty)$: 
\[ H(x,t)\leq \widehat{H}(t) =
\max_{x\in M}H(x,0)\,e^{(2n+\frac{1}{n})t} < +\infty.\]

This shows that the harmonic mean curvature flow exists on
$[0,+\infty)$.
} \qed

In the rest of this section, we do not assume the ambient
manifold $N^{n+1}$ is a hyperbolic manifold.

%
\beginproposition \label{second order estimate of F}
{Assume $N^{n+1}$ is a smooth $n+1\geq 3$ dimensional manifold
which is convex at infinity, the maximal existence time of the 
harmonic mean curvature flow $\phi:M\times [0,T)\to N$ is $T=
+\infty$, and as $t\to +\infty,$ $M_t=\phi(M,t)$ converges to a
smooth $n$ dimensional submanifold $M_{\infty}$ of $N$ in the
$C^{\infty}$-topology; then \[\max_{x\in M,\, t\in
[0,+\infty)}\{|F(x,t)|, |\nabla F(x,t)|, |\nabla^2 F(x,t)|\}\leq
C,\] where $C$ is a constant depending on $M_0$, $N^{n+1}$  and
$M_{\infty}$.
}
\end{proposition}

\pf
{ Straightforward from the assumptions.
} \qed

%
\begin{proposition}\label{the integral estimate of F} 
{Assume $N$ and $M_t\rightarrow M_{\infty}$ are as in the hypotheses of
Proposition \ref{second order estimate of F}. Then
\[\lim_{t\rightarrow \infty}\int_{M_t}F^2\,  d\mu_t=0.\]
}
\end{proposition}

\pf
{ By Theorem $1.1$ in \cite{EHR}, we have the formula
$\frac{\partial}{\partial t}(\int_{M_t}d\mu_t)=
-\int_{M_t}FHd\mu_t$. Because
$\int_{M_t}d\mu_t\rightarrow \mu(M_{\infty})$ as
$t\rightarrow \infty$, we could find an $\epsilon$-dense set
$\{t_k\}_{k=1}^{\infty}$ for any positive constant $\epsilon>0$
such that \[\lim_{k\rightarrow \infty}t_k=\infty\] and
\[\lim_{k\rightarrow \infty}\int_{M_{t_k}}FH\, d\mu_{t_k}=0.\]
Then using the inequality $H\geq n^2 F$, we get
$\lim_{k\rightarrow \infty}\int_{M_{t_k}}F^2\, d\mu_{t_k}=0$.\\

Now to get our conclusion we only need to show
$\frac{\partial}{\partial t}\int_{M_t}F^2d\mu_t$ is uniformly
bounded. First, we know from Proposition \ref{second order
estimate of F} that $|F|$,
$|\nabla F|$ and $|\nabla^2 F|$ are uniformly bounded. So we have
\[\frac{\partial}{\partial t} \Big(\int_{M_t} F^2 \, d\mu_t \Big)
=\int 2FF_t+F^2(-FH)\, d\mu_t\]

\[=\int 2F\Big(\mathscr L(F)+ \sum_{i=1}^n F \Big(\frac{\partial f}{\partial
\lambda_i} \Big) (\lambda_i^2+\mathscr R_{ii}) \Big)- F^3H\, d\mu_t\]
(where we use equation (\ref{formula}))
\[=\int 2nF^4+ 2F^4 \Big(\sum_{i=1}^n\lambda_i^{-2}\mathscr R_{ii} \Big)
+2F\mathscr L(F)-F^3H\, d\mu_t\]
\[\leq \int 2F^4K_2 \Big(\sum_{i=1}^n\lambda_i^{-2} \Big)\, d\mu_t+
\int 2F\mathscr L(F)\, d\mu_t\]
\[\leq C\int F^2\, d\mu_t+2\int F\mathscr L(F)\, d\mu_t,\]
where the first inequality uses the following facts:\\
$M_t$ is always contained in some compact set of $N^{n+1}$, since
$N^{n+1}$ is convex at infinity, so its sectional curvature is
bounded above by some constant $K_2$;  and
$HF^{-1}=(\sum_{i=1}^n\lambda_i)(\sum_{i=1}^n\lambda_i^{-1})\geq
n^2\geq 2n$.\\

Next, since we know the volume of $M_t$ is always non-increasing and 
$|F|$ is uniformly bounded, we get
\[C\int_{M_t} F^2\, \, d\mu_t\leq C_1,\]
where $C_1$ is some constant depending only on $M_0$, $N$ and
$M_{\infty}$.\\
Since $|\nabla^2 F|$ is uniformly bounded, we get
\[2\int F\mathscr L(F)\, d\mu_t\leq
2n^2\int F|\nabla^2 F|\, d\mu_t\leq C_2,\]
where $C_2$ is some constant depending on $M_0$, $N$ and
$M_{\infty}$.\\
By all the above we get
\[\frac{\partial}{\partial t}\Big(\int F^2\, d\mu_t\Big)\leq
C_3,\]
where $C_3$ is another constant depending on $M_0$, $N$ and
$M_{\infty}$.\\
Therefore
\[\lim_{t\rightarrow \infty}\int_{M_t}F^2\, d\mu_t=0.\]
} \qed

%
\begin{corollary}
{Assume $N$ and $M_t\rightarrow M_{\infty}$ are as assumed for
Proposition \ref{second order estimate of F}. Then we have
\[\lim_{t\rightarrow \infty}\Big(\max_{x\in M}F(x,t)\Big)=0.\]
}
\end{corollary}

\pf
{ By Proposition \ref{the integral estimate of F}, we have
\[0= \lim_{t\rightarrow \infty} \int_{M_t} F^2 \, d\mu_t=
\int_{M_\infty}\lim_{t\rightarrow\infty} F^2(x,t) \,
d\mu_{\infty},
\]

\noindent
so the corollary follows.
} \qed

By the above results, assume $N$ and $M_t\rightarrow M_{\infty}$
are as in the hypotheses of Proposition \ref{second order estimate
of F}, we know that $F\equiv 0$ on the limit surface $M_{\infty}$,
if $M_{\infty}$ is the smooth limit of the harmonic mean curvature
flow, which implies that $\det \mathscr W=0$ on $M_{\infty}$.\\

%
%
\section{Classification of harmonic mean curvature flow on
surfaces}

In this section, we consider harmonic mean curvature flow for
$n=2$, where $M^2$ is an orientable surface, $N^3$ is a hyperbolic
manifold, and the harmonic mean
$f(\lambda)= \frac{\lambda_1\lambda_2}{\lambda_1+ \lambda_2}$.
As before, we assume that $\phi_0(M)$ is locally strictly convex.

In the following we always assume $F(x,0)< \frac12$, i.e.
$\lambda_1^{-1}+\lambda_2^{-1}> 2$, which will guarantee the
harmonic mean curvature flow exists forever by Lemma
\ref{hyperbolic infinite time}.  Note that, for example,
$f(\lambda_1, \lambda_2) <\frac12$ for the examples of Theorems
\ref{dimred} and \ref{nodimred}, and that the horospheres have
$f(\lambda_1, \lambda_2) \equiv\frac12$.\\

We define $C_0=2\pi \chi(M_0)=\int_{M_t}(K-1) \, d\mu_t$,
where the second equation is true for any $M_t$ because
of the Gauss-Bonnet theorem, where $\chi(M_0)$
is the Euler number of $M_0$;
$K(x,t)=\lambda_1(x,t) \lambda_2(x,t)$,
$\lambda_1(x,t)$ and $\lambda_2(x,t)$ are the
principal curvatures at the point $x$ on $M_t$ in the ambient
hyperbolic
manifold $N^3$; and the Gauss equation, which implies the Gauss
curvature $= K-1$.\\

First, define $V(t)=\int_{M_t}1 \, d\mu_t$, the area of $M_t$.
Then using the formula
\[\frac{\partial}{\partial t}d\mu_t= -FH d\mu_t\]
we get
\[\frac{d}{dt}V(t)
=\int_{M_t}\frac{\partial}{\partial t}\, d\mu_t
=\int_{M_t}(-FH)\, d\mu_t=\int_{M_t}(-K)\, d\mu_t\]
\[=-\int_{M_t}(K-1)\, d\mu_t -\int_{M_t}1 \, d\mu_t
=-C_0-V(t).\]
Solving the above ODE, we get
\[V(t)=(V(0)+ C_0)e^{-t}- C_0.\]
This shows that the area of $M_t$ is determined by its genus and
the area $V(0)$ of the initial surface $M_0$.

There are three cases: $C_0<0$, $C_0=0$, $C_0>0$, corresponding to
the surfaces with genus $g>1$ (Case I) , $g=1$ (Case II) and
$g=0$ (Case III). \\

(I). Let us first consider the case $C_0= 2\pi \chi(M_0)< 0$. In
this case, we have
\[\lim_{t \rightarrow \infty}V(t)=-C_0> 0\]
which means the limit surface has non-zero volume. We
{\bf conjecture} that in a hyperbolic manifold $N^3$, the
limit surface will be the totally geodesic surface, if there is
one
in the homotopy class of $M_0$.  This behavior is seen in Theorem
\ref{nodimred}.\\

(II). When $C_0= 2\pi \chi(M_0) =0$, we have
\[\lim_{t\rightarrow \infty}V(t)=-C_0= 0\]
which means the limit surface has zero volume. In fact we could
prove the following:

%
\begin{proposition}
{If $N^3$ is a hyperbolic manifold, $F(x,0)< \frac12$ for all
$x\in M$ and the genus of $M$ $=0$, then
\[\lim_{t\rightarrow \infty} (\max_{x\in M_t}H(x,t)) = +\infty. \]
}
\end{proposition}

\pf
{ Because $\int_{M_t}(K-1)\, d\mu_t = C_0=0$, we have
$\max_{x\in M_t}K(x,t)\geq 1$.
We also have
$\lim_{t \rightarrow \infty} (\max_{x\in M_t}F(x,t))= 0$,
using the assumption $F(x,0)< \frac12$, by
Lemma \ref{hyperbolic infinite time}. Then for any $x\in M_t$,
$t>0$, we have the following: \[K(x,t)=H(x,t)F(x,t)\leq
F(x,t)(\max_{x\in M_t} H(x,t)).\] Taking the maximum on the both
sides of the above inequality, we have

\[1\leq \max_{x\in M_t}K(x,t)\leq
(\max_{x\in M_t}F(x,t))(\max_{x\in M_t} H(x,t)).\]
So
\[\max_{x\in M_t} H(x,t) \geq \frac{1}{\max_{x\in M_t}F(x,t)}.\]
Taking the limit on both sides, we get
\[\lim_{t \rightarrow \infty}(\max_{x\in M_t} H(x,t)) \geq
\frac{1}{\lim_{t \rightarrow \infty}(\max_{x\in M_t}F(x,t))}=
+\infty.\]
} \qed

The above proposition means that there exists at least one blowup
point on the limit set;  the example of Theorem \ref{dimred} blows
up at every point.\\

(III) Finally, when $C_0= 2\pi \chi(M_0) >0$, we have an
interesting geometric result. In this case, because
\[V(t)=(V(0)+ C_0)e^{-t}- C_0,\] there exists some
$T_0$, $0<T_0< +\infty$,
such that $V(T_0)=0$. That means the harmonic mean curvature flow
stops in finite time. But we already proved that the flow will
exist forever if $F< \frac{1}{2}$.  So under the assumption
$F< \frac12$, this surface will not exist.

%
\begin{remark}
{Observe that the nonexistence of the initial surfaces in Case
(III) above may also be proven by lifting the simply-connected
surface $M_0$ to the universal cover $H^3$ of $N^3$ and applying
the comparison principle with shrinking spheres centered at a
point: the sphere of radius $r$ has 
$F=\frac12 \coth r > \frac12.$}
\end{remark}

%
%
\section{General geometric flows}

In this section we give examples for a general geometric flow
\eqref{eq:1} in a hyperbolic manifold $N^{n+1}$ which will
exist forever or for a computable finite time, and converge to a 
given totally geodesic submanifold $P^k$ of any codimension.  
In this section, we always assume the existence of a totally 
geodesic submanifold $P^k$ in $N^{n+1}$.

Firstly, by similar methods to those of section $2$ and section 
$3$, we may prove a theorem for general dimensions and codimensions:

%
\begin{theorem}\label{P^k}
{ Assume $P^k$ is a compact totally geodesic submanifold of the
hyperbolic manifold $N^{n+1}$, where $1\leq k\leq n$. Let $M$ be 
diffeomorphic to the unit sphere bundle of the normal bundle 
$\perp P$\ when $k<n$ ; we choose $M$ to be one of the two 
connected components of the unit sphere bundle of the normal
bundle $\perp P$ when $k=n$.  Then we have a flow by
harmonic mean curvature $\phi_t: M\rightarrow N$ such that as
$t\rightarrow +\infty$, $\phi_t(M)\rightarrow P$.  }
\end{theorem}

\pf
{ We only sketch the proof.  We find the second
fundamental form matrix of $\psi_r(M)$ with respect to a basis of
curvature directions is the following:

\begin{equation*}
\mathscr{W} = \left(
\begin{array}{cc}
I_k \tanh r & 0_{k\times (n-k)}  \\
0_{(n-k) \times k} & I_{n-k} \coth r
\end{array} \right)
\end{equation*}

Then we find
\begin{equation}\label{dr/dt=HMCF}
 \frac{\partial r}{\partial t}= -F= -\frac{\tanh r}{k+
(n-k)(\tanh r)^2}
\end{equation}
Solving this ODE, we get
\[(\sinh r(t))^k(\cosh r(t))^{n-k}= C e^{-t}\]
where $C= (\sinh r_0)^k(\cosh r_0)^{n-k}$ is a fixed positive
constant.  This shows that $\phi_t := \psi_{r(t)}$ is a solution
of harmonic mean curvature flow.

Note that $r(t) \rightarrow 0$ as $t \rightarrow +\infty$.
} 
\qed

Now let $M^n$ be diffeomorphic to (one connected
component of) the unit sphere normal bundle of $P^k$ in $N^{n+1}$,
and let $\psi_r:M \to N$ define the hypersurface at distance $r>0$
from $P^k$.  We consider flow by an arbitrary symmetric function
of the normal curvatures:

%
\begin{theorem} \label{last}
{For the symmetric function $f(\lambda_1 , \cdots ,\lambda_n)$,
define
$$h(r) = f(\tanh r, \cdots ,\coth r),$$ 
where $\tanh r$ is repeated $k$ times and $\coth r$ is repeated
$n-k$ times.  Choose $r_0>0$ and define
\[T_0 = \int_0^{r_0} \frac{1}{h(r)} dr, \quad
0 < T_0 \leq +\infty. \]
Then we may construct a flow 
\begin{equation}\label{flow}
{\frac{\partial}{\partial t} \phi(\cdot , t)= 
f(\lambda(\mathscr{W}(x,t)))\,\vec{v}(x,t)
}
\end{equation}
with initial condition $\phi(\cdot , 0) = \psi_{r_0}$, which exists for
time $0\leq t\leq T_0\leq\infty$, and $\phi(\cdot , t)$ converges to the
totally geodesic $k$-dimensional submanifold $P^k$ as $t \to T_0$.
}
\end{theorem}

\pf
{The hypersurface defined by $\phi(\cdot,t) := \psi_{r(t)}$ flows
by \eqref{flow} if
\begin{equation}\label{dr/dt=h}
\frac{\partial r}{\partial t}= -F(x,t) \equiv -h(r)
\end{equation}
\[\Longrightarrow \int_{r(0)}^{r(T_0)} \frac{1}{h(r)} dr=
\int_{0}^{T_0} -1\, dt\]
\[\Longrightarrow T_0=\int_0^{r_0} \frac{1}{h(r)} dr.\]

The conclusion now follows from the proof of Theorem \ref{P^k},
replacing equation \eqref{dr/dt=HMCF} with equation
\eqref{dr/dt=h}.
} \qed

%
\begin{remark}
{Note that the flow (\ref{flow}) is parabolic if 
$\frac{\partial f}{\partial \lambda_i}> 0$ $(1\leq i \leq n)$; 
parabolic for backwards time if 
$\frac{\partial f}{\partial \lambda_i}< 0$ $(1\leq i \leq n)$; 
and is a first-order PDE if $f$ is constant.  }
\end{remark}

The following corollary is a generalization of both mean curvature
flow ($m=1$, $\ell=0$) and of harmonic mean curvature
flow ($m=n$, $\ell=n-1$).

%
\begin{corollary}\label{S/S}
{Assume $P^k$ is a compact totally geodesic submanifold of
$N^{n+1}$, where $1\leq k\leq n$. Let $M$ be diffeomorphic to 
the unit sphere bundle of the normal bundle $\perp P$ when $k<n$; 
$M$ is one of the two components of the unit sphere bundle of 
$\perp P$ when $k=n$.

For integers $0\leq m$, $\ell\leq n$, let $S_m$, $S_\ell$ be the 
elementary symmetric functions of degree $m$, $\ell$ respectively, of
the principal curvatures $\lambda_1, \dots, \lambda_n$ of $M_t$.
We have a flow by curvature function
$$F(x,t)=
\frac{S_m(\lambda_1, \cdots , \lambda_n)}
{S_\ell(\lambda_1, \cdots , \lambda_n)},$$
for time $0\leq t<\infty$,
such that $\phi(t): M\rightarrow N,$ and $\phi_t(M)\rightarrow P$
as $t\rightarrow +\infty$; assuming that the integers $m$, $\ell$
satisfy $|m-(n-k)| < |\ell-(n-k)|$.
}
\end{corollary}

%
\begin{remark}
{Theorem \ref{last} also may be applied to prove a partial
converse of Corollary \ref{S/S}:   assuming $P^k$ and $N^{n+1}$
are as in Corollary \ref{S/S}, if the opposite condition
$|m- (n-k)|\geq |\ell- (n-k)|$ holds, then the same
construction
yields a flow of hypersurfaces by the curvature function
$F = \frac{S_m}{S_\ell}$ which converges to the totally geodesic
submanifold $P^k$ in {\bf finite} time $T_0$.  }
\end{remark}

\pf
{ In the following, we fix an arbitrary positive constant $r(0)= r_0$.
Firstly we have

\begin{equation*}
S_m =
\sum^n_{\substack{p+q=m \\ 0\leq p\leq k \\ 0\leq q\leq n-k}}
C_{k}^{p} (\tanh r)^p \, C_{n-k}^{q} (\coth r)^q =
\sum C_{k}^{p}C_{n-k}^{q} (\coth r)^{q- p},
\end{equation*}
where $C_k^p$ is the combinatorial coefficient $\frac{k!}{p!
(k-p)!}$.

Since $\coth r\geq 1$, it is easy to see
\begin{equation*}
S_m \sim \left\{
\begin{array}{rl}
(\coth r)^{m} & \text{if } m\leq n-k\\
(\coth r)^{2(n-k)- m} & \text{if } m> n-k
\end{array} \right.
\end{equation*}
where the notation $S_m \sim (\coth r)^j$ means that there exist
positive constants $C_1$, $C_2$ such that
$C_1(\coth r)^j \leq S_m \leq C_2(\coth r)^j$. Here $C_1$ and $C_2$
will depend only on $m$, $n$, $k$, $\ell$ and $r_0$.

Similarly, we have

\begin{equation*}
S_\ell \sim \left\{
\begin{array}{rl}
(\coth r)^{\ell} & \text{if } \ell\leq n-k\\
(\coth r)^{2(n-k)- \ell} & \text{if } \ell> n-k.
\end{array} \right.
\end{equation*}

Therefore

\begin{equation*}
F=\frac{S_m}{S_\ell} \sim \left\{
\begin{array}{rl}
(\coth r)^{m- \ell} & \text{if } m, \ell\leq n-k\\
(\coth r)^{\ell- m} & \text{if } m, \ell> n-k\\
(\coth r)^{2(n-k)- m- \ell} & \text{if } \ell\leq n-k < m \\
(\coth r)^{m+ \ell- 2(n-k)} & \text{if } m\leq n-k <  \ell.
\end{array} \right.
\end{equation*}

By Theorem \ref{last}, we obtain that the flow exists forever if
and only if the power of $\coth r$ is negative in the asymptotic 
estimate for $F$ above.  That is, if and only if $m$ and $\ell$
satisfy one of the following conditions:

\begin{equation*}
\left\{
\begin{array}{rl}
m< \ell & \text{if } m, \ell\leq n-k\\
\ell< m & \text{if } m, \ell> n-k\\
2(n-k)< m+ \ell & \text{if } \ell \leq n-k < m \\
m+ \ell< 2(n-k) & \text{if } m\leq n-k < \ell.
\end{array} \right.
\end{equation*}

It is straightforward to see the above inequalities are equivalent
to the inequality $|m-(n-k)| < |\ell-(n-k)|$, which is our conclusion.
} \qed

%
\begin{remark}
{In particular, the case $k=n,\ m=1,\ \ell=0$ is the first example
we are aware of in the literature of a locally {\bf convex} compact 
hypersurface
flowing by mean curvature and converging smoothly to a submanifold in
infinite time. And the case $k= n-1,\ m=0,\ \ell=1$ gives an
example of (backwards parabolic)
inverse mean curvature flow existing forever and
converging to a totally geodesic hypersurface. After reversing
time to obtain parabolicity, this example of $-\frac{1}{H}$ flow
is properly divergent as $t \to \infty$.}
\end{remark}

\end{document}